\documentclass[12pt]{article} 

\usepackage{amssymb}
\usepackage{amsmath}
\usepackage{a4wide}
\usepackage{latexsym}
\usepackage{mathrsfs}
\usepackage[v2,tips]{xy}
\usepackage[T1]{fontenc}
\usepackage[latin1]{inputenc}
\usepackage{theorem}

\renewcommand{\mathcal}{\mathscr}

\theoremstyle{change}

\begingroup \makeatletter

\gdef\th@changenop{\normalfont\itshape
  \def\@begintheorem##1##2{%
        \item[\hskip\labelsep \theorem@headerfont
          \mathversion{bold}##2\mathversion{normal}]}%
\def\@opargbegintheorem##1##2##3{%
\item[\hskip\labelsep \theorem@headerfont
  \mathversion{bold}##2\ ##3\mathversion{normal}]}}

\endgroup

\newtheorem{theorem}{Theorem.}[section]
\newtheorem{proposition}[theorem]{Proposition.}
\newtheorem{lemma}[theorem]{Lemma.}

\newtheorem{corollary}[theorem]{Corollary.}
{\theorembodyfont{\rmfamily}\newtheorem{definition}[theorem]{Definition.}}
{\theorembodyfont{\rmfamily}}
{\theorembodyfont{\rmfamily}\newtheorem{remark}[theorem]{Remark.}}
{\theorembodyfont{\rmfamily}}
{\theorembodyfont{\rmfamily}}
{\theorembodyfont{\rmfamily}}
{\theoremstyle{changenop}
  \theorembodyfont{\rmfamily}}

\newcommand{\comment}[1]{{}}

\numberwithin{equation}{section}

\newcounter{smallromans}

\newenvironment{romanenumerate}
{\begin{list}{{\normalfont\textrm{(\roman{smallromans})}}}
  {\usecounter{smallromans}\setlength{\itemindent}{0cm}
   \setlength{\leftmargin}{5.5ex}\setlength{\labelwidth}{5.5ex}
   \setlength{\topsep}{0.75\parsep}\setlength{\partopsep}{0ex}
   \setlength{\itemsep}{0ex}}}
{\end{list}}

\newcommand{\romanref}[1]{{\normalfont\textrm{(\ref{#1})}}}

\newcounter{smallromansdash}

\newcounter{smallromansprime}

\newcounter{smallarabics}

\newcounter{smallalphs}

\newcommand{\allop}{\ensuremath{\mathcal{B}}}
\newcommand{\weaklycompactop}{\ensuremath{\mathcal{W}}}
\newcommand{\compactop}{\ensuremath{\mathcal{K}}}
\newcommand{\ssop}{\ensuremath{\mathcal{S}}}
\newcommand{\inessop}{\ensuremath{\mathcal{E}}}
\newcommand{\ccop}{\ensuremath{\mathcal{V}}}
\newcommand{\opidealg}{\ensuremath{\mathcal{G}}}
\newcommand{\closedopidealg}{\ensuremath{\overline{\mathcal{G}}}}

\newcommand{\indexset}[1]{\ensuremath{\Gamma_{{#1}}}}

\newcommand{\uc}{\operatorname{uc}}
\newcommand{\card}{\operatorname{card}}

\newcommand{\spa}{\operatorname{span}}
\newcommand{\closedspa}{\ensuremath{\overline{\spa}\,}}

\newcommand{\N}{\ensuremath{\mathbb{N}}}

\newcommand{\R}{\ensuremath{\mathbb{R}}}
\newcommand{\C}{\ensuremath{\mathbb{C}}}
\newcommand{\K}{\ensuremath{\mathbb{K}}}
\newcommand{\J}{\ensuremath{\mathbb{J}}}

\renewcommand{\phi}{\ensuremath{\varphi}}
\renewcommand{\epsilon}{\ensuremath{\varepsilon}}

\renewcommand{\leq}{\ensuremath{\leqslant}}
\renewcommand{\le}{\ensuremath{\leqslant}}
\renewcommand{\geq}{\ensuremath{\geqslant}}
\renewcommand{\ge}{\ensuremath{\geqslant}}

\newcommand{\smashw}[2][l]{{\text{\makebox[0pt][#1]{$#2$}}}}

\newcommand{\beginpf}{\smallskip%
\noindent\textsl{Proof. }}
\newcommand{\eopf}{\hfill $\Box$}

\author{Alistair Bird, Graham Jameson and Niels Jakob
  Laustsen}

\title{The Giesy--James theorem for general index~$p$, with an
  application to operator ideals on the $p^{\text{th}}$ James space}

\begin{document}

\maketitle

\begin{abstract} \noindent A theorem of Giesy and James states that
  $c_0$ is finitely representable in James' quasi-reflexive Banach
  space~$J_2$. We extend this theorem to the $p^{\text{th}}$
  quasi-reflexive James space~$J_p$ for each $p\in(1,\infty)$.  As an
  application, we obtain a new closed ideal of operators on~$J_p$,
  namely the closure of the set of operators that factor through the
  complemented subspace
  $(\ell_\infty^1\oplus\ell_\infty^2\oplus\cdots\oplus
  \ell_\infty^n\oplus\cdots)_{\ell_p}$ of~$J_p$.\\[2mm]
  \noindent 2010 \emph{Mathematics Subject Classification:}  Primary
  46B45, 47L20; Secondary 46B07, 46H10, 47L10.\\[2mm]
  \noindent
  \emph{Key words and phrases:} quasi-reflexive Banach space, James
  space, finite repre\-sen\-tability of~$c_0$, closed operator ideal.
\end{abstract}

\section{Introduction}
As outlined in the abstract, we shall prove that $c_0$ is finitely
representable in the $p^{\text{th}}$ quasi-reflexive James space~$J_p$
for each $p\in(1,\infty)$ and then show how this result gives rise to
a new closed ideal of operators on~$J_p$. In order to make these
statements precise, let us introduce some notation and terminology.

We denote by $\N_0$ and $\N$ the sets of non-negative and positive
integers, respectively.  Following Giesy and James~\cite{gj},
we index sequences by~$\N_0$ and write $x(n)$ for the $n^{\text{th}}$
element of the sequence~$x$, where $n\in\N_0$.  For a non-empty subset
$A$ of~$\N_0$, we write $A = \{n_1<n_2<\cdots<n_k\}$ (or \mbox{$A =
  \{n_1<n_2<\cdots\}$} if $A$ is infinite) to indicate that
$\{n_1,n_2,\ldots,n_k\}$ is the increasing ordering of~$A$.

Let $\K=\R$ or $\K=\C$ be the scalar field, and let
$p\in(1,\infty)$. For a scalar sequence $x$ and a finite subset~$A=
\{n_1<n_2<\cdots<n_{k+1}\}$ of~$\N_0$ of cardinality at least two, we
define
\[ \nu_p(x,A) = \Bigl(\sum_{j=1}^k \bigl| x(n_j) -
x(n_{j+1})\bigr|^p\Bigr)^{\frac{1}{p}}; \] for convenience, we let
$\nu_p(x,A) = 0$ whenever $A\subseteq\N_0$ is empty or a singleton.
Then $\nu_p(\,\cdot\,,A)$ is a seminorm on the vector
space~$\K^{\N_0}$ of all scalar sequences, and
\begin{align*}
  \|x\|_{J_p} := \sup&\bigl\{\nu_p(x,A):
  A\subseteq\N_0,\, \card A<\infty\bigr\}\\
  =\sup&\biggl\{ \Bigl(\sum_{j=1}^k \bigl| x(n_j) -
  x(n_{j+1})\bigr|^p\Bigr)^{\frac{1}{p}} : k\in\N,\,n_1,\ldots,n_{k+1}\in
  \N_0,\, n_1 <\cdots < n_{k+1}\biggr\}
\end{align*}
defines a complete norm on the subspace $J_p := \bigl\{ x\in c_0 : \|
x\|_{J_p} < \infty\bigr\}$, which we call the \emph{$p^{\text{th}}$
  James space}. The sequence $(e_m)_{m=0}^\infty$,
where~$e_m\in\K^{\N_0}$ is given by 
\[ e_m(n) = \begin{cases} 1 &\text{if}\ m=n\\ 0
  &\text{otherwise} \end{cases}\qquad (n\in\N_0), \] forms a shrinking
Schauder basis for~$J_p$.  More importantly, $J_p$ is
\emph{quasi-reflexive} in the sense that the canonical image of~$J_p$
in its bidual has codimension one. This result, as well as the
definition of~$J_p$, is due to James~\cite{james} in the case $p=2$;
Edelstein and Mityagin~\cite{em} appear to have been the first to
observe that it carries over to arbitrary $p\in(1,\infty)$.

A Banach space $X$ is \emph{finitely representable} in a Banach
space~$Y$ if, for each finite-dimensional subspace~$F$ of~$X$ and each
$\epsilon>0$, there is an operator $T\colon F\to Y$ such that
\begin{equation}\label{eq_finiterepr}
  (1-\epsilon)\|x\|_X\le\|Tx\|_Y\le (1+\epsilon)\|x\|_X\qquad (x\in
  F). \end{equation}
We shall in fact only consider finite representability of~$c_0$, in
which case 
it suffices to establish~\eqref{eq_finiterepr} for the
finite-dimensional subspaces $F = \ell_\infty^n$, where
\mbox{$n\in\N$}.  
Although not required, let us mention the Maurey--Pisier theorem
that $c_0$ is finitely representable in a Banach space~$Y$ if
and only if $Y$ fails to have finite cotype (\emph{e.g.},
see~\cite[Theorem~14.1]{djt}). This result shows  in particular 
that finite representability of~$c_0$ is an isomorphic invariant,
despite the obvious dependence on the 
choice of norm in~\eqref{eq_finiterepr}.

Giesy and James~\cite{gj} proved that $c_0$ is finitely
representable in~$J_2$. Our first main result, to be proved in
Section~\ref{proofofgengjthm}, extends this result to arbitrary
$p\in(1,\infty)$.

\begin{theorem}\label{giesyjamesthm}
  For each $p\in(1,\infty)$, $c_0$ is finitely representable in~$J_p$.
\end{theorem}

To explain how this result leads to a new closed ideal of operators
on~$J_p$, we require some more notation.  For $p\in[1,\infty)$ and a
family $(X_j)_{j\in\J}$ of Banach spaces, we write
$\bigl(\bigoplus_{j\in\J} X_j\bigr)_p$ for the direct sum of the
$X_j$'s in the sense of~$\ell_p$; that is,
\[ \Bigl(\bigoplus_{j\in\J} X_j\Bigr)_p = \Bigl\{ (x_j) : x_j\in X_j\
(j\in\J)\ \ \text{and}\ \ \sum_{j\in\J} \|x_j\|^p <\infty\Bigr\}. \]
We shall only apply this notation in two cases, namely
\begin{equation}\label{defnlpsumoflinftyns}
  G_p := \Bigl(\bigoplus_{n\in\N} \ell_\infty^n\Bigr)_p\qquad
  \text{and}\qquad J_p^{(\infty)} := \Bigl(\bigoplus_{n\in\N_0}
  J_p^{(n)}\Bigr)_p, 
\end{equation}
where $J_p^{(n)}$ denotes the subspace of~$J_p$ spanned by the
first~$n+1$ basis vectors $e_0,e_1,\ldots,e_n$.  \newpage%
\noindent
Our interest in these spaces stems from the two facts that (i)~$J_p$
contains a complemented subspace isomorphic to~$J_p^{(\infty)}$; and
(ii)~Theorem~\ref{giesyjamesthm} implies that $J_p^{(\infty)}$
contains a complemented subspace isomorphic to~$G_p$ (for $p=2$, this
has already been observed by Casazza, Lin and
Lohman~\cite[Theorem~13(i)]{cll} using the original Giesy--James
theorem), and this subspace gives rise to a new closed ideal of
operators on~$J_p$, as we shall now outline.

For Banach spaces $X$ and $Y$, let
\[ \opidealg_Y(X) = \bigl\{ ST : T\in\allop(X,Y),\,
S\in\allop(Y,X)\bigr\} \] be the set of operators on~$X$ which factor
through~$Y$. This defines a two-sided algebraic ideal of the Banach
algebra~$\allop(X)$ of bounded operators on~$X$, provided that~$Y$
contains a complemented subspace isomorphic to~$Y\oplus Y$ (which will
always be the case in this paper), and hence its norm-closure, denoted
by~$\closedopidealg_Y(X)$, is a closed ideal of~$\allop(X)$.

Edelstein and Mityagin~\cite{em} made the easy, but fundamental,
observation that the quasi-reflexivity of~$J_p$ for $p\in(1,\infty)$
implies that the ideal~$\weaklycompactop(J_p)$ of weakly compact
operators has codimension one in~$\allop(J_p)$, hence is a maximal
ideal.  Loy and Willis~\cite[Open Problems~2.8]{lw} formally raised
the problem of determining the structure of the lattice of closed
ideals of~$\allop(J_2)$, having themselves proved that
$\compactop(J_2)\subsetneq\closedopidealg_{\ell_2}(J_2)\subsetneq
\weaklycompactop(J_2)$ and $\ssop(J_2) =
\inessop(J_2)\not\supseteq\closedopidealg_{\ell_2}(J_2)$, where
$\ssop(J_2)$ and $\inessop(J_2)$ denote the ideals of strictly
singular and inessential operators, respectively
(see~\cite[Theorem~2.7]{lw} and the text preceding it). Saksman and
Tylli~\cite[Remark~3.9]{st} improved the latter result by showing that
$\compactop(J_2) = \ssop(J_2)$, while the third
author~\cite{lau1,lau2} generalized these results to arbitrary
$p\in(1,\infty)$ and, more importantly, complemented them by showing
that the lattice of closed ideals in $\allop(J_p)$ has the following
structure:
\begin{equation*}
  \spreaddiagramrows{-2.7ex} \xymatrix{ 
    \allop(J_p)\ar@{-}[d]\\
    \smashw[r]{\weaklycompactop(J_p)}=%
    \opidealg_{J_p^{(\infty)}}(J_p)%
    =\smashw[l]{\closedopidealg_{J_p^{(\infty)}}(J_p)}%
    \ar@{.}[d]\\
    \closedopidealg_{\ell_p}(J_p)\ar@{-}[d]\\
    \smashw[r]{
      \compactop(J_p)=\ssop(J_p)}=%
    \smashw[l]{\inessop(J_p)=\ccop(J_p)}\ar@{-}[d]\\
    \{0\}\smashw[l]{,}} 
\end{equation*}
where 
$\ccop(J_p)$ is the ideal of completely continuous operators, the
vertical lines indicate proper set-theoretic inclusion, and further
closed ideals may be found only at the dotted line. In particular,
$\weaklycompactop(J_p)$ is the unique maximal ideal of~$\allop(J_p)$.

The second main result of this paper, which we shall prove in
Section~\ref{section_pfofThmnewoperatoridealonJp}, states that
$\allop(J_p)$ contains at least one other closed ideal than those
listed above.

\begin{theorem}\label{newoperatoridealonJp}
  For each $p\in(1,\infty)$, the operator ideal
  $\closedopidealg_{G_p}(J_p)$ lies strictly between
  $\closedopidealg_{\ell_p}(J_p)$ and $\weaklycompactop(J_p)$,
  where~$G_p$ is the Banach space given
  by~\eqref{defnlpsumoflinftyns}.

  Hence the lattice of closed ideals in~$\allop(J_p)$ has at least
  six distinct elements, namely
  \[ \{0\}\subsetneq
  \compactop(J_p)\subsetneq\closedopidealg_{\ell_p}(J_p)\subsetneq
  \closedopidealg_{G_p}(J_p)\subsetneq \weaklycompactop(J_p)\subsetneq
  \allop(J_p).
\] 
\end{theorem}

\section{Proof of Theorem~\ref{giesyjamesthm}}%
\label{proofofgengjthm}

Throughout this section, we fix a number $p\in (1,\infty)$.  Our aim
is to prove Theorem~\ref{giesyjamesthm} by modifying the proof of
Giesy and James \cite{gj}.  The general scheme of the proof is the
same, but at several points, identities that are simple in the case $p
= 2$ have to be replaced with estimations applying to other $p$.  We
follow their notation as far as possible.  We show that there is a
near-isometric embedding of $\ell_\infty^K$ for each $K\in\N$ in the
real case.  It then follows easily, by standard techniques, that there
is at least an isomorphic embedding in the complex case.

Spiky vectors play a central role in the proof. As in \cite[p.\
65]{gj}, let
\[ z_{2k} = \frac{1}{(2k)^{1/p}} \sum_{j=1}^k e_{2j-1}\in J_p\qquad
(k\in\N), \] so that $z_{2k}$ is a unit vector with spikes in its
initial $k$ odd coordinates.

The other key ingredient is the ``stretch'' operator $T_n\colon J_p\to
J_p$ which, for $n\in\N$ and $x\in J_p$, is given by $(T_nx)(kn) =
x(k)$ whenever $k\in\N_0$ and by linear interpolation in between these
points.  One can easily check that~$T_n$ is linear and isometric.

We use the notation $[j,k]$ for the set of integers $n$ such that $j
\leq n \leq k$.

By an inductive process, we construct, for each $K\in\N$, a set of $K$ 
stretched spiky vectors with the parameters chosen suitably, and show
that these vectors are equivalent to the usual basis of
$\ell_\infty^K$.  The inductive step is captured by the following
lemma, corresponding to \cite[Lemma~1]{gj}.

\begin{lemma}\label{mainlemma}
  Let $m\in\N$ and $\gamma,\epsilon\in(0,\infty)$.  Suppose that $x$
  is an element of $J_p$ supported on the integer interval $[0, 2m-1]$
  and satisfying
  \begin{equation}\label{condx2} 
    \max_{0\le j< 2m}\bigl|x(j) - x(j+1)\bigr|^p \le
    \frac{\gamma}{2m}\qquad \text{and}\qquad 
    \| x \|_{J_p}^p - \nu_p\bigl(x, [0,2m]\bigr)^p\le \epsilon.
  \end{equation}
  For some even $n$, let $w = T_n x + \gamma^{1/p}z_{2mn}$.  Then $w$
  is supported on the integer interval $[0, 2mn -1]$ and satisfies
  \begin{equation}\label{condw1} 
    \max_{0 \leq j < 2mn} | w(j) - w(j+1) |^p\le \frac{\gamma }{2mn} 
    \Bigl( 1 + \frac{1}{n^{1 - 1/p}} \Bigr)^p 
  \end{equation}
  and
  \begin{equation}\label{condw2} 
    \| w\|_{J_p}^p - \nu_p \bigl( w, [0,2mn] \bigr)^p\le 2\epsilon +
    \gamma \phi(m,n),
  \end{equation}
  where $\phi(m,n) \to 0$ as $n \to \infty$ with $m$ fixed.
\end{lemma}

We show next how Theorem~\ref{giesyjamesthm} follows, and then return
to the proof of Lemma \ref{mainlemma}.
\bigskip

\noindent 
\textsl{Proof of Theorem~\ref{giesyjamesthm}.}  With $\epsilon > 0$
and $K \in \N$ given, we construct vectors $x_1,\ldots,x_K\in J_p$
with $\| x_i \|_{J_p} \geq 1$ for $1\le i\le K$ such that $\|
\sum_{i=1}^K \delta_i x_i \|_{J_p} \leq 1 + 2\epsilon$ for all choices
of \mbox{$\delta_1,\ldots,\delta_K\in\{-1, 1\}$}.  We then deduce
equivalence with the usual basis of $\ell_\infty ^K$ as follows.  By
convexity, we have $\bigl\| \sum_{i=1}^K \lambda_i
x_i\bigr\|_{J_p}\leq 1 + 2\epsilon$ for all real $\lambda_i$ with
$|\lambda_i| \leq 1$.  Suppose that $\max_{1 \leq i \leq K}
|\lambda_i| = |\lambda_j| = 1$.  Then $\bigl\|\sum_{i=1}^K \lambda_i
x_i - 2\lambda_j x_j\bigr\|_{J_p} \leq 1+ 2\epsilon$ (the coefficient
of $x_j$ has been changed to $-\lambda_j$), so
\[ \Bigl\| \sum_{i=1}^K \lambda_i x_i \Bigr\|_{J_p} \geq \| 2\lambda_j
x_j \|_{J_p} - \Bigl\| \sum_{i=1}^K \lambda_i x_i - 2\lambda_j x_j
\Bigr\|_{J_p} \geq 2 - (1 + 2\epsilon ) = 1 - 2\epsilon. \]

Let $\epsilon_k = \epsilon /3^{K - k}$.  At stage $k$, we will define
$n_k\in\N$, $\gamma_k\in\R$ and \mbox{$x_1^{(k)},\ldots,x_k^{(k)}\in
  J_p$} such that the following properties hold.  Firstly,
$x_1^{(k)},\ldots,x_k^{(k)}$ are supported on the integer interval
\mbox{$[0, 2m_k-1]$}, where $m_k := n_1 n_2 \ldots n_k$, and $\|
x_i^{(k)}\|_{J_p} \geq 1$ for $1\le i\le k$.  Secondly,
\[ 1\le \gamma_k \leq 1 + \frac{\epsilon k}{K}. \] Thirdly, for all
choices of $\delta_1,\ldots,\delta_k \in \{-1, 1\}$ and with $y_\delta
^{(k)} := \sum_{i=1}^k \delta_i x_i^{(k)}$, we have
\begin{equation}\label{ind2}
  \max_{0\leq j< 2m_k} \bigl|y_\delta ^{(k)} (j) - y_\delta ^{(k)}(j+1)\bigr|^p 
  \leq \frac{\gamma_k}{2m_k}
\end{equation}
and
\begin{equation}\label{ind3}
  \| y_\delta^{(k)} \|_{J_p}^p -  \nu_p \bigl( y_\delta^{(k)}, [0,
  2m_k] \bigr)^p\le \epsilon_k. 
\end{equation}
By \eqref{ind2} and \eqref{ind3}, we then obtain $\| y_\delta ^{(k)}
\|_{J_p}^p \leq \gamma_k + \epsilon_k \leq 1 + 2\epsilon\leq (1 +
2\epsilon)^p$, from which the desired conclusion follows.

To start, take $x_1^{(1)} = z_2$ and $n_1 = \gamma_1 = 1$.  Suppose
now that stage $k-1$ has been completed.  For a certain even
integer~$n_k$ to be chosen, define
\[ x_i^{(k)} = T_{n_k} (x_i^{(k-1)})\qquad (1\leq i \leq k-1)
\quad\qquad\mbox{and}\quad\qquad x_k^{(k)} = \gamma_{k-1}^{{1}/{p}}
z_{2m_k}. \] Let $\delta_1,\ldots,\delta_k\in\{-1,1\}$ be given.  We
may assume that $\delta_k = 1$.  Apply Lemma \ref{mainlemma} with $x =
y_\delta^{(k-1)}$, $m = m_{k-1}$, $n = n_k$, $\epsilon =
\epsilon_{k-1}$ and $\gamma = \gamma_{k-1}$.  Then
\[ w = T_{n_k} (y_\delta ^{(k-1)}) + \gamma_{k-1}^{{1}/{p}}
z_{2m_k} = y_\delta ^{(k)}, \] hence~\eqref{condw1} implies
that~\eqref{ind2} is satisfied with
\[ \gamma_k = \gamma_{k-1} \biggl(1 + \frac{1}{n_k^{1- 1/p}}
\biggr)^p. \] We choose $n_k$ large enough to ensure that $\gamma_k
\leq 1 + \epsilon k/K$.  By \eqref{condw2},
\[ \| y_\delta^{(k)} \|_{J_p}^p - \nu_p \bigl(y_\delta^{(k)}, [0,
2m_k] \bigr)^p \leq 2 \epsilon_{k-1} + \gamma_{k-1} \phi(m_{k-1},
n_k). \] Since $\epsilon_k = 3\epsilon_{k-1}$, to ensure \eqref{ind3},
we choose $n_k$ also to satisfy $ \gamma_{k-1} \phi(m_{k-1}, n_k) \leq
\epsilon_{k-1}$.  \eopf
 
\begin{remark} Because of the dependence of $\phi(m,n)$ on $m$, it is
  not possible to take $n_k$ equal to the same value $n$ for each $k$,
  as in \cite{gj} for the case $p = 2$.  We shall actually see later
  that $\phi(m,n)$ only depends on $m$ when $p > 2$.
\end{remark}

\noindent
\textsl{Outline of proof of Lemma \ref{mainlemma}.}  Write $y = T_n x$
and $z = \gamma^{1/p}z_{2mn}$, so that $w = y+z$.  Clear\-ly, $y$ and
$z$ are both supported on the integer interval $[0, 2mn -1]$.  Also,
from the definitions, we have \mbox{$|z(j) - z(j+1)| = (\gamma
  /2mn)^{1/p}$} and
\[ \bigl|(T_nx)(j)-(T_nx)(j+1)\bigr|\le
\frac{1}{n}\Bigl(\frac{\gamma}{2m}\Bigr)^{\frac{1}{p}}=
n^{\frac{1}{p}-1}\Bigl(\frac{\gamma}{2mn}\Bigr)^{\frac{1}{p}}\qquad
(0\le j < 2mn), \] from which \eqref{condw1} follows.

The bulk of the work is the proof of \eqref{condw2}.  Since $w$ is
supported on the integer interval $[0, 2mn-1]$, we can find a set $A =
\{ a_1< a_2< \cdots< a_{k+1}\}$, with $a_1 = 0$ and $a_{k+1} = 2mn$,
such that $\| w\|_{J_p} =\nu_p(w,A)$.  The aim is to show that the
whole interval acts as a reasonable substitute for this set~$A$.  This
will be accomplished by four steps, summarized as follows:
\begin{align} 
  \nu_p(w,A)^p &\leq \nu_p(y,A)^p + \nu_p(z,A)^p +
  \rho_1\label{step1}\\ 
  &\leq \nu_p\bigl(y,A\cup ([0,2mn]\cap n\N_0)\bigr)^p +
  \nu_p\bigl(z,A\cup([0,2mn]\cap n\N_0)\bigr)^p + \rho_1 +
  \rho_2\label{step2}\\ 
  &\leq \nu_p\bigl(y,[0,2mn]\bigr)^p + \nu_p\bigl(z,[0,2mn]\bigr)^p +
  \rho_1 + \rho_2\label{step3}\\ 
  &\leq \nu_p\bigl(w, [0,2mn]\bigr)^p + \rho_1 + \rho_2 +
  \rho_3,\label{step4}
\end{align} 
where $\rho_1$, $\rho_2$ and~$\rho_3$ are error terms which will
emerge from the proofs.  Step 1 moves from $w = y+z$ to $y$ and $z$
separately, and step 4 reverses this.  Working with $y$ and $z$
separately, step 2 adjoins multiples of $n$ to $A$, and step 3 adjoins
all intervening integers.  Because of the concepts involved, we
present these four steps in the order 1, 4, 3, 2.  \eopf

\begin{lemma}\label{ineq1}
  Suppose that $a, \, b > 0$.  Then $(a+b)^p - a^p - b^p \leq 2^p
  (a^{p-1}b + ab^{p-1})$.
\end{lemma}
\beginpf
With no loss of generality, we may assume that $a \geq b$.  Writing
$b/a = t$, we see that the stated inequality is equivalent to $(1+t)^p
- 1 - t^p \leq 2^p(t + t^{p-1})$ for $0 < t \leq 1$.  For such~$t$,
since the function $t \mapsto (1+t)^p$ is convex and $t = (1-t).0 +
t.1$, we have $(1+t)^p \leq (1-t).1 + 2^p t$, hence $(1+t)^p - 1 \leq
(2^p -1)t$, which of course implies the required inequality.  \eopf

\begin{remark}
  The estimation in Lemma~\ref{ineq1} is quite adequate for our
  purposes. In fact, the best constant on the right-hand side of the
  inequality is $p$ for $2 \leq p \leq 3$, and $2^{p-1}-1$
  otherwise~\cite{jameson}.
\end{remark}

\noindent
\textsl{Step 1: Proof of} \eqref{step1}, with 
\begin{equation}\label{eq_rho1} \rho_1 = \begin{cases} \displaystyle{2^p\gamma
      \Bigl(\frac{1}{n^{p - 2 + 1/p}} + \frac{1}{n^{1 - 1/p}}\Bigr)} &
    \text{for}\ 1<p\le
    2,\\[2.5ex]
    \displaystyle{2^p\gamma \Bigl(\frac{(2m)^{p-2}}{n^{1/p}} + \frac{1}{n^{1-
          1/p}} \Bigr)} & \text{for}\ p>2. \end{cases} \end{equation}
Write $\ell_i = a_{i+1} - a_i$, so that $\sum_{i=1}^k \ell_i = 2mn$.
Then we have, by definition,
\[ \bigl|y(a_i)-y(a_{i+1})\bigr| \leq \frac{\ell_i}{n}
\Bigl(\frac{\gamma}{2m}\Bigr)^{\frac{1}{p}} \qquad \text{and}\qquad
\bigl|z(a_i)-z(a_{i+1})\bigr|\leq
\Bigl(\frac{\gamma}{2mn}\Bigr)^{\frac{1}{p}}\qquad (1\le i\le k). \]
Lemma \ref{ineq1} implies that $\nu_p(y+z, A)^p - \nu_p(y, A)^p -
\nu_p(z, A)^p \leq 2^p s$, where
\begin{align*}
  s\,\smashw{:}&= \sum_{i=1}^k
  \Bigl(\bigl|y(a_i)-y(a_{i+1})\bigr|^{p-1} \bigl|z(a_i)-
  z(a_{i+1})\bigr| + \bigl|y(a_i)-y(a_{i+1})\bigr|
  \,\bigl|z(a_i)-z(a_{i+1})\bigr|^{p-1}
  \Bigr)\\
  &\leq \sum_{i=1}^k \biggl(\Bigl(\frac{\ell_i}{n} \Bigr)^{p-1}
  \Bigl(\frac{\gamma }{2m}\Bigr)^{1 - \frac{1}{p}}
  \Bigl(\frac{\gamma}{2mn}\Bigr)^{\frac{1}{p}} +
  \frac{\ell_i}{n}\Bigl(\frac{\gamma }{2m}\Bigr)^{\frac{1}{p}}
  \Bigl(\frac{\gamma }{2mn}\Bigr)^{1 - \frac{1}{p}}\biggr)\\
   &= \gamma\biggl(\frac{\sum_{i=1}^n \ell_i^{p-1}}{2mn^{p-1+1/p}} +
  \frac{1}{n^{1- 1/p}}\biggr),
\end{align*}
since $\sum_{i=1}^k \ell_i = 2mn$.  
For $1 < p \leq 2$, we have $\ell_i^{p-1} \leq \ell_i$, hence
\[ s \leq \gamma \Bigl(\frac{1}{n^{p - 2 + 1/p}} + \frac{1}{n^{1 -
    1/p}}\Bigr), \] whereas for $p > 2$, $\sum_{i=1}^k \ell_i^{p-1}
\leq \bigl(\sum_{i=1}^k \ell_i\bigr)^{p-1} = (2mn)^{p-1}$, so that 
\[ s \leq \gamma \Bigl(\frac{(2m)^{p-2}}{n^{1/p}} + \frac{1}{n^{1-
    1/p}} \Bigr). \] Multiplying these upper bounds on~$s$ by~$2^p$,
we conclude that~\eqref{step1} is satisfied with $\rho_1$ given
by~\eqref{eq_rho1}.  \eopf \bigskip

\noindent
\textsl{Step 4: Proof of} \eqref{step4}, with $\rho_3 =
\gamma/n^{p-1}$.  
Letting 
\[ s_k = \sum_{j=kn}^{(k+1)n-1}\Bigl(\bigl|w(j) - w(j+1)\bigr|^p -
\bigl|y(j) - y(j+1)\bigr|^p - \bigl|z(j) - z(j+1)\bigr|^p
\Bigr)
, \] we can write
\[ \nu_p\bigl(w, [0,2mn]\bigr)^p - \nu_p\bigl(y,[0,2mn]\bigr)^p -
\nu_p\bigl(z,[0,2mn]\bigr)^p = \sum_{k=0}^{2m-1}s_k. \]  Our claim is
that this quantity is at least $-\gamma/n^{p-1}$.

To verify this, fix integers $k\in [0,2m-1]$ and $j\in [kn, (k+1)n
-1]$. Then $z(j) - z(j+1)$ is alter\-nate\-ly~$\pm c$, where $c :=
(\gamma /2mn)^{1/p}$, while $y(j) - y(j+1) = \frac{1}{n}\bigl(x(k) -
x(k+1)\bigr)$, and by assumption \mbox{$d_k := \frac{1}{n}\bigl| x(k)
  - x(k+1)\bigr|\leq \frac{1}{n} (\gamma/2m)^{1/p}$}.
Since $c > d_k$, we see that \mbox{$\bigl|w(j) - w(j+1)\bigr|$} is
alternately $c + d_k$ and $c - d_k$. Hence, as $n$ is even,
\begin{equation}\label{eqsk}
  s_k = \frac{n}{2}\bigl((c + d_k)^p + (c - d_k)^p - 2d_k^p -
  2c^p\bigr). \end{equation}  
By convexity of the function $t \mapsto t^p$, we have $(c + d_k)^p +
(c - d_k)^p \geq 2c^p$. Therefore
\mbox{$s_k \geq -nd_k^p \geq -\gamma/2mn^{p-1}$}, so 
$\sum_{k=0}^{2m-1}s_k \geq -\gamma/n^{p-1}$, as required. 
\eopf

\begin{remark}
  Equation~\eqref{eqsk} shows that~$s_k=0$ for $p = 2$, and in fact
  one can prove that $s_k\ge 0$ whenever $p \geq 2$, thus rendering
  the error term~$\rho_3$ superfluous for such~$p$.
\end{remark}

We now come to Step~3, which is really the heart of the method, and it
is the one where it is essential to work with
$\nu_p(\,\cdot\,,\,\cdot\,)^p$ rather than
$\nu_p(\,\cdot\,,\,\cdot\,)$ itself.  We shall adjoin all integers to
the set $A\cup([0,2mn]\cap n\N_0)$.  This has the effect of reducing
$\nu_p(y,\,\cdot\,)^p$, but the reduction is more than offset by an
increase in $\nu_p(z, \,\cdot\,)^p$.

\begin{lemma}\label{ineq2}
  Suppose that $t \geq 1$. Then $t^p - t \leq (t-1)(t+1)^{p-1}$.
\end{lemma}
\beginpf
For $1 < p \leq 2$, we have $t^{p-1}\le t$, hence $t^p - t \leq t^p -
t^{p-1} = (t-1)t^{p-1}$, which is stronger than the stated inequality.
For $p \geq 2$, we use the convexity of the function $t \mapsto
t^{p-1}$.  Since
\[ t = \frac{t-1}{t} (t+1) + \frac{1}{t}.1,\] we have
\[ t^{p-1} \leq \frac{t-1}{t} (t+1)^{p-1} + \frac{1}{t},\] which is
again stronger than the stated inequality.  \eopf \bigskip

\noindent
\textsl{Step~3: Proof of~\eqref{step3}.}  Let $B = A\cup([0,2mn]\cap
n\N_0)$. Our aim is to prove that
\[ \delta := \nu_p\bigl(y,[0,2mn]\bigr)^p +
\nu_p\bigl(z,[0,2mn]\bigr)^p - \nu_p(y,B)^p - \nu_p(z,B)^p \] is
non-negative.  Writing $B = \{b_1 < b_2 <\cdots < b_{h+1}\}$, we have
\mbox{$\delta = \sum_{j=1}^h\bigl(\Delta_j(y) + \Delta_j(z)\bigr)$},
where
\[ \Delta_j(y) = \biggl(\sum_{i=b_j}^{b_{j+1}-1}\bigl|y(i) -
y(i+1)\bigr|^p\biggr) - \bigl|y(b_j) - y(b_{j+1})\bigr|^p \] and
$\Delta_j (z)$ is defined similarly. Hence it suffices to prove that
$\Delta_j(y) + \Delta_j(z)\ge 0$ for each integer $j\in[1,h]$.

By the choice of~$B$, $b_j$ and $b_{j+1}$ both belong to an interval
of the form~$[kn,(k+1)n]$ for some $k\in\N$. As in the proof of
Step~4, this implies that
\[ \bigl|y(i) - y(i+1)\bigr| = d_k\qquad (b_j\le i < b_{j+1})\qquad
\text{and}\qquad \bigl|y(b_j) - y(b_{j+1})\bigr| = \ell_j d_k, \]
where \mbox{$d_k := \frac{1}{n}\bigl| x(k) - x(k+1)\bigr|\leq
  \frac{1}{n} (\gamma/2m)^{1/p}$} and $\ell_j := b_{j+1} - b_j$, so
\mbox{$\Delta_j(y) = (\ell_j - \ell_j^p) d_k^p$}.  Meanwhile,
$\bigl|z(i) - z(i+1)\bigr| = c$ for each~$i
$, where $c := (\gamma /2mn)^{1/p}$, and $\bigl|z(b_j) -
z(b_{j+1})\bigr|$ equals 0 if $\ell_j$ is even and $c$ if $\ell_j$ is
odd, thus in both cases $\Delta_j(z) \geq (\ell_j - 1)c^p$.

Now if $\ell_j \leq n-1$, we find
\[ \Delta_j(y) + \Delta_j(z) \geq (\ell_j - \ell_j^p)d_k^p + (\ell_j -
1)c^p \geq \bigl((\ell_j - \ell_j^p) + (\ell_j -1) n^{p-1}\bigr)d_k^p \]
because $c^p\geq n^{p-1}d_k^p$.  Since $n \geq \ell_j + 1$, Lemma
\ref{ineq2} gives $(\ell_j -1) n^{p-1} \geq \ell_j^p - \ell_j$, hence
$\Delta_j(y) + \Delta_j(z) \geq 0$, as required.  Otherwise $\ell_j =
n$, which is assumed even, so that $\Delta_j(z) = n c^p$, and
\begin{xxalignat}{3}
  && \Delta_j(y) + \Delta_j(z) &= (n - n^p)d_k^p + nc^p \geq (n - n^p
  + n^p)d_k^p = n d_k^p > 0. &&\Box
\end{xxalignat} 

Finally, we reach Step~2 where multiples of~$n$ are adjoined to the
set~$A$. We require two lemmas, the first of which describes the
effect on $\nu_p(\,\cdot\,, A)^p$ of substituting new end points
in~$A$, while the second considers the effect of filling in gaps in~$A$.

\begin{lemma}\label{endpoints}
  Consider integers $\ell\geq 3$ and $0\le c \leq b_1 < b_2 < \cdots
  < b_\ell\le c'$, let $B = \{ b_1,b_2,\ldots,b_\ell\}$ and $C =
  \{ c, b_2, \ldots,b_{\ell-1},c' \}$, and suppose that $v\in
  J_p$ satisfies
  \begin{equation}\label{endpointsEq1} v(c) \leq v(b_1) \leq v(b_j)
    \leq v(b_\ell) \leq v(c')\quad \text{or}\quad v(c) \geq
    v(b_1) \geq v(b_j) \geq v(b_\ell) \geq v(c')\ \mbox{} \end{equation} 
  for $1 < j < \ell$.  Then
  \[ \nu_p\bigl(v, \{b_1, b_\ell\}\bigr)^p - \nu_p(v,B)^p \leq \nu_p
  \bigl(v,\{ c, c'\}\bigr)^p - \nu_p(v,C)^p. \]
\end{lemma}

\beginpf
We consider only the case where the first set of inequalities
in~\eqref{endpointsEq1} is satisfied; the other case is similar.  We
replace the end points of~$B$ one at a time.  Let $D = \{ c, b_2,
\ldots ,b_{\ell-1}, b_\ell \}$.  In the sum under consideration, $r :=
v(b_2) - v(b_1)$ is replaced with $s := v(b_2) -v(c)$, and both are
non-negative, so $\nu_p(v, D)^p - \nu_p(v, B)^p = s^p -
r^p$. Differen\-tia\-tion shows that the function $t \mapsto (s+t)^p -
(r+t)^p$, is increasing on~$[0,\infty)$ because $s\ge r$, and hence
$s^p - r^p\le (s+t)^p - (r+t)^p$ for each $t\ge 0$. Taking $t :=
v(b_\ell) - v(b_2)$, we obtain $s+t = v(b_\ell) - v(c)$ and $r+t =
v(b_\ell) - v(b_1)$, so 
\[ \nu_p(v, D)^p - \nu_p(v, B)^p \leq \nu_p\bigl(v,\{ c,
b_\ell\}\bigr)^p - \nu_p\bigl(v, \{b_1, b_\ell\}\bigr)^p. \] A similar
argument with $r := v(b_\ell) - v(b_{\ell-1})$, $s := v(c') -
v(b_{\ell-1})$ and \mbox{$t := v(b_{\ell-1}) - v(c)$} shows that
  \[ \nu_p(v, C)^p - \nu_p(v, D)^p\leq \nu_p\bigl(v,\{ c,
  c'\}\bigr)^p - \nu_p\bigl(v, \{c, b_\ell\}\bigr)^p. \] Adding
  these two inequalities, we conclude that
\[ \nu_p(v, C)^p - \nu_p(v, B)^p \leq \nu_p\bigl(v,\{ c,
c'\}\bigr)^p - \nu_p\bigl(v, \{b_1, b_\ell\}\bigr)^p, \] from
which our statement follows.  \eopf

\begin{lemma}\label{fillgaps}
  Let $\ell\in\N$, and suppose that $C_1,\ldots,C_\ell$ and
  $D_1,\ldots,D_\ell$ are finite subsets of $\N_0$ with \mbox{$\min
    C_j = \min D_j =: m_j$} and $\max C_j = \max D_j =: m'_j$, where
  $m'_j \leq m_{j+1}$ for each~$j$.  Suppose further that
  $E_1,\ldots,E_{\ell-1}$ are finite subsets of~$\N_0$ such that $\min
  E_j = m'_j$ and $\max E_j = m_{j+1}$ for each $j$ (so $E_j$ is
  between $C_j\cup D_j$ and $C_{j+1}\cup D_{j+1}$), and let $E_\ell=\{
  m'_\ell \}$.  Then
  \[ \sum_{j=1}^\ell \bigl(\nu_p(v,D_j)^p - \nu_p(v,C_j)^p \bigr) =
  \nu_p \Bigl(v, \bigcup_{j=1}^\ell (D_j \cup E_j)\Bigl)^p -
  \nu_p\Bigl(v, \bigcup_{j=1}^\ell (C_j \cup E_j)\Bigr)^p\qquad (v\in
  J_p). \]
\end{lemma}
\beginpf
Clearly, we have \[ \nu_p \Bigl(v, \bigcup_{j=1}^\ell (C_j \cup
E_j)\Bigl)^p = \sum_{j=1}^\ell \nu_p (v, C_j)^p + \sum_{j=1}^\ell
\nu_p(v, E_j)^p, \] which together with the corresponding formula for
$\nu_p\bigl(v, \bigcup_{j=1}^\ell (C_j \cup E_j)\bigr)^p$ gives the
result. \mbox{} \eopf \bigskip

\noindent
\textsl{Step 2: Proof of}~\eqref{step2}, with $\rho_2 =
2\epsilon$. Let $N = [0,2mn]\cap n\N_0$.  The effect on $z$ of
adjoining elements to the set \mbox{$A = \{ a_1<\cdots < a_{k+1}\}$}
is easily seen.  Let $\ell_i = a_{i+1} - a_i$ for $1\le i\le k$.  As
in the proof of Step~3 above, $\bigl|z(a_i) - z(a_{i+1})\bigr|$ is
$c:=(\gamma /2mn)^{1/p}$ if $\ell_i$ is odd, and~$0$ if $\ell_i$ is
even.  If $\ell_i$ is odd and new points are inserted between $a_i$
and $a_{i+1}$, then at least one of the new intervals, say $[b_j,
b_{j+1}]$, has odd length, so $\bigl|z(b_j) - z(b_{j+1})\bigr| = c$.
Hence \[ \nu_p(z, A) \leq \nu_p (z, A\cup N). \]

We shall now prove the corresponding inequality for $y$, just with an
error term added on the right-hand side.  Recall that $a_1 = 0$ and
$a_{k+1} = 2mn$.  Note that if there is some $b\in N$ such that $a_i <
b < a_{i+1}$ for some~$i$ and either \mbox{$y(b)<\min\bigl\{y(a_i),
  y(a_{i+1})\bigr\}$} or \mbox{$y(b)>\max\bigl\{y(a_i),
  y(a_{i+1})\bigr\}$}, then \mbox{$\bigl|y(a_i) - y(b)\bigr|^p +
  \bigl|y(b) - y(a_{i+1})\bigr|^p > \bigl|y(a_i) -
  y(a_{i+1})\bigr|^p$}.  Hence we may adjoin any such points~$b$ to
the set $A$, thereby increasing $\nu_p(y, A)$ without changing $A\cup
N$; we still use the notation \mbox{$A= \{a_1<\cdots < a_{k+1}\}$} for
the augmented set.

Let the intervals $[a_i, a_{i+1}]$\ $(1\le i\le k)$ that contain at
least one multiple of $n$ be relabelled $[b_j, b_j']$\ $(1\le j\le h)$
and ordered increasingly; that is, $b_1< b'_1\le b_2 <b'_2\le\cdots\le
b_h<b'_h$.  Note that $b_1 = a_1 = 0$ and $b'_h = a_{k+1} = 2mn$, and
that $b'_j$ may or may not be equal to~$b_{j+1}$ for $j\le h-1$.
Then, with $B_j := \bigl([b_j,b'_j]\cap n\N\bigr)\cup\{b_j, b'_j\}$
for $1\le j\le h$, we have
\begin{equation}\label{step2eq1}
  \nu_p(y, A)^p - \nu_p(y, A \cup N)^p = \sum_{j=1}^h
  \bigl(\nu_p\bigl(y,\{ b_j, b'_j \}\bigr)^p - 
  \nu_p(y, B_j)^p \bigr). \end{equation} 
Let $c_j = \max\bigl([0,b_j]\cap n\N_0\bigr)$ and
$c'_j = \min\bigl([b_j',\infty)\cap n\N_0\bigr)$, and let $C_j
= [c_j, c'_j] \cap n\N_0$. Then $c_1 = 0$ and $c'_h = 2mn$, and  Lemma
\ref{endpoints} implies that 
\begin{equation}\label{step2eq2}
 \nu_p \bigl(y, \{ b_j, b'_j \}\bigr)^p - \nu_p (y, B_j)^p \leq
\nu_p \bigl(y, \{ c_j, c'_j \}\bigr)^p - \nu_p (y, C_j)^p\qquad (1\le
j\le h). \end{equation}
(Note that the augmentation of the set~$A$ carried out in
the previous paragraph ensures that~$y$ satisfies the
hypothesis~\eqref{endpointsEq1}.)

We now seek to invoke Lemma~\ref{fillgaps} with the
sets~$\{c_j,c_j'\}$ playing the role of the $D_j$'s. To do so, we
require some more notation. Let $c_0=c_0' = 0$, $C_0 = \{0\}$,
$c_{h+1} = c'_{h+1} = 2mn$ and $C_{h+1} = \{2mn\}$.  Then clearly
$\min\{c_j,c_j'\} = \min C_j = c_j$ and $\max\{c_j,c_j'\} = \max C_j =
c_j'$ for each integer~$j\in[0,h+1]$, but $c_j'\le c_{j+1}$ need not
be satisfied for each $j\le h$. It is, however, true that $c_j'\le
c_{j+2}$ for each $j\le h-1$ (because the interval
$[b_{j+1},b_{j+1}']$ contains a multiple of~$n$). Hence, taking $E_j =
[c_j',c_{j+2}]\cap n\N_0$ for $0\le j\le h-1$ and letting $E_h =
E_{h+1} = \{2mn\}$, we can apply Lemma~\ref{fillgaps} for even and odd
indices~$j$ separately.  We observe that $C_j\cup E_j = [c_j, c_{j+2}]
\cap n\N_0$ for $0\le j\le h-1$, so $\bigcup_{j\in \indexset{r}}
(C_j\cup E_j) = N$ for $r\in\{0,1\}$, where $\indexset{0}$ and
$\indexset{1}$ denote the sets of even and odd integers in~$[0,h+1]$,
respectively.  Thus Lemma~\ref{fillgaps} gives
\begin{equation}\label{oddsum} 
  \sum_{j\in \indexset{r}} \bigl(\nu_p\bigl(y, \{ c_j, c'_j\}\bigr)^p -
  \nu_p(y, C_j)^p \bigr) = \nu_p \Bigl(y, \bigcup_{j\in
    \indexset{r}}(\{c_j,c_j'\}\cup E_j)\Bigr)^p - \nu_p (y, N)^p. 
\end{equation}
Since $y=T_n x$ and $N = [0,2mn]\cap n\N_0$, we have $\nu_p (y, N)^p =
\nu_p\bigl(x, [0,2m]\bigr)^p \geq \| x\|_{J_p}^p - \epsilon$
by~\eqref{condx2}, while $\nu_p\bigl(y, \bigcup_{j\in
  \indexset{r}}(\{c_j,c_j'\}\cup E_j)\bigr)^p \leq \| y\|_{J_p}^p = \|
x\|_{J_p}^p$.  Hence the sum in \eqref{oddsum} is no greater than
$\epsilon$, so adding the two cases ($r=0$ and $r=1$) and
using~\eqref{step2eq1} and~\eqref{step2eq2}, we conclude that
\begin{xxalignat}{3}  &&\nu_p(y, A)^p - \nu_p(y, A \cup N)^p &\leq
\sum_{j=1}^h\bigl(\nu_p\bigl(y, \{ c_j, c'_j\}\bigr)^p - \nu_p(y,
C_j)^p \bigr)\\ &&&= 
\sum_{j=0}^{h+1}\bigl(\nu_p\bigl(y, \{ c_j, c'_j\}\bigr)^p - \nu_p(y,
C_j)^p \bigr)
\le2\epsilon. &&\Box \end{xxalignat}

\noindent
\textsl{Completion of the proof of Lemma \ref{mainlemma}.}  With the
four steps completed, it is clear that Lemma 2.1 holds with
\[ \phi(m,n) = 2^p \Bigl(\psi(m,n) + \frac{1}{n^{1 - 1/p}}\Bigr) +
\frac{1}{n^{p-1}}, \] where
\[ \psi(m,n) = \begin{cases} \displaystyle{\frac{1}{n^{p -2 + 1/p}}} &
  \text{for}\ 1 < p \leq
  2\\[2.5ex]
  \displaystyle{\frac{(2m)^{p-2}}{n^{1/p}}} & \text{for}\ p > 2.
\end{cases} \] Note that $m$ does not appear in the case $p \leq 2$,
and that $p - 2 + 1/p > 0$, so in both cases $\phi(m,n) \to 0$ as $n \to
\infty $ with $m$ fixed.  \eopf

\section{Proof of Theorem~\ref{newoperatoridealonJp}}%
\label{section_pfofThmnewoperatoridealonJp}

We begin with an elementary observation which is tailored to reduce
Theorem~\ref{newoperatoridealonJp} to the statement that the Banach
spaces $G_p$ and~$J_p^{(\infty)}$ given by~\eqref{defnlpsumoflinftyns}
are non-isomorphic. A closely related result can be found
in~\cite[Proposition~5.3.8]{pie}. 
\begin{lemma}\label{lemmaopidealsandcomplementedsubspaces}
  Let $X$, $Y$ and $Z$ be Banach spaces satisfying:
  \begin{romanenumerate}
  \item $X$ contains a complemented subspace isomorphic to~$Y;$
  \item $Y$ contains a complemented subspace isomorphic to~$Z;$
  \item $Y\cong Y\oplus Y$ and $Z\cong Z\oplus Z$.
  \end{romanenumerate}
  Then $\closedopidealg_Z(X)\subseteq\closedopidealg_Y(X)$, with
  equality if and only if $Z\cong Y$.
\end{lemma}

\beginpf The inclusion
$\closedopidealg_Z(X)\subseteq\closedopidealg_Y(X)$ is clear, as is
the equality of these two sets in the case where $Z\cong Y$.

Conversely, suppose that $\closedopidealg_Z(X) =
\closedopidealg_Y(X)$, and let~$P$ be a projection on~$X$ with
$P(X)\cong Y$. Clearly $P$ factors through~$Y$, so $P$ belongs
to~$\closedopidealg_Z(X)$ by the assumption. It then follows from
standard results that~$Z$ contains a complemented subspace isomorphic
to~$Y$ (\emph{e.g.}, see~\cite[Propostion~3.4 and Lemma~3.6(ii)]{lau1}
for details), and therefore $Y$ and~$Z$ are isomorphic by the
Pe{\l}czy{\a'n}ski decomposition method. \eopf \bigskip

We shall next record the facts required to invoke
Lemma~\ref{lemmaopidealsandcomplementedsubspaces} in the proof of
Theorem~\ref{newoperatoridealonJp}. 
\begin{lemma}\label{factsaboutcomplsubspacesoJp} For each
  $p\in(1,\infty)$,
  \begin{romanenumerate}
  \item\label{factsaboutcomplsubspacesoJp1} $G_p$ contains a
    complemented subspace isomorphic to~$\ell_p;$
  \item\label{factsaboutcomplsubspacesoJp2} $J_p^{(\infty)}$ contains
    a complemented subspace isomorphic to~$G_p;$
  \item\label{factsaboutcomplsubspacesoJp3} $J_p$ contains a
    complemented subspace isomorphic to~$J_p^{(\infty)};$
  \item\label{factsaboutcomplsubspacesoJp4} $\ell_p\cong
    \ell_p\oplus\ell_p$, $G_p\cong G_p\oplus G_p$ and 
    $J_p^{(\infty)}\cong J_p^{(\infty)}\oplus J_p^{(\infty)}$.
  \end{romanenumerate}
\end{lemma}

\beginpf All but one of these results are well known. The exception
is~\romanref{factsaboutcomplsubspacesoJp2} which, however, follows
from Theorem~\ref{giesyjamesthm} in exactly the same way as the
corresponding result for $p=2$ is deduced from the original
Giesy--James theorem in~\cite[Theorem~13(i)]{cll}. 

References for the other statements are as follows;
\romanref{factsaboutcomplsubspacesoJp1} and the first part
of~\romanref{factsaboutcomplsubspacesoJp4} are obvious,
while~\romanref{factsaboutcomplsubspacesoJp3} and the remaining two
parts of~\romanref{factsaboutcomplsubspacesoJp4} follow
from~\cite[Lemmas~5 and~6]{em}. (A key condition appears to be missing
in the statement of \cite[Lemma~5]{em}, though, namely that the
sequence denoted by~$\nu$ is unbounded.) \eopf

\begin{remark}\label{remark_bddbelowoperators} Let $X$ and $Y$ be
  Banach spaces.  An operator $T\colon X\to Y$ is \emph{bounded below}
  by~$\epsilon>0$ if $\|Tx\|_Y\ge \epsilon\|x\|_X$ for each $x\in X$. In
  this case~$T$ is an isomorphism onto its image, and the inverse
  operator has norm at most~$\epsilon^{-1}$, so in particular the
  Banach--Mazur distance~$d_{\text{BM}}$ between the domain~$X$ and
  the image~$T(X)$ of~$T$ satisfies
  \[ d_{\text{BM}}\bigl(X,T(X)\bigr)\le\frac{\|T\|}{\epsilon}. \]

  Now suppose that $X$ is a closed subspace of~$Y$ and that $T\colon
  X\to Y$ is linear and satisfies
  \begin{equation*}
    \| x - Tx\|\le \eta\|x\|\qquad (x\in X)
  \end{equation*}
  for some $\eta\in(0,1)$. Then we have $(1-\eta)\|x\|\le
  \|Tx\|\le(1+\eta)\|x\|$ for each $x\in X$, so by the previous
  paragraph $T$ is an isomorphism onto its image, and
  \[ d_{\text{BM}}\bigl(X,T(X)\bigr)\le \frac{1+\eta}{1-\eta}. \]
\end{remark}

\begin{definition}
Let~$F$ be a finite-dimensional Banach
  space. 
  The \emph{unconditional basis constant} of a basis
  $\mathbf{b}=\{b_1,\ldots,b_n\}$ for~$F$  is given by
  \[ K_{\mathbf{b}} := \sup\biggl\{\Bigl\|\sum_{j=1}^n \alpha_j\beta_j
  b_j\Bigr\| : \alpha_j,\beta_j\in\K,\, |\alpha_j|\le 1\
  (j=1,\ldots,n),\, \Bigl\|\sum_{j=1}^n \beta_j b_j\Bigr\|\le
  1\biggr\}. \] The infimum of the unconditional basis constants of
  all possible bases for~$F$ is the \emph{unconditional constant}
  of~$F$; we denote it by~$\uc(F)$.
\end{definition}

It is easy to verify that, for Banach spaces~$E$ and $F$ of the same
finite dimension, we have
\begin{equation}\label{remarklusteq1} \uc(E)\le
  d_{\text{BM}}(E,F)\uc(F). 
\end{equation}

\begin{definition}\label{defnlust} (Dubinsky,
Pe{\l}czy{\a'n}ski and Rosenthal~\cite[Definition~3.1]{dpr}.)
  Let $C\in[1,\infty)$. A Banach space~$X$ has \emph{local
    unconditional structure} (or \emph{l.u.st.}\ for short) with
  \emph{constant} at most~$C$ if each finite-dimensional subspace
  of~$X$ is contained in some larger finite-dimensional subspace~$F$
  of~$X$ with $\uc(F)\le C$.
\end{definition}

A Banach space with an unconditional basis has l.u.st. This applies in
particular to~$G_p$. On the other hand, Johnson and
Tzafriri~\cite[Corollary~2]{jt} have shown that no quasi-reflexive
Banach space has l.u.st. We shall use this result to prove that
$J_p^{(\infty)}$ does not have l.u.st.

We begin with a generalization of the above-mentioned fact that every
Banach space with an unconditional basis has l.u.st. This result is
probably well-known to specialists, but as we have been unable to locate a
reference, we include a proof.
\begin{lemma}\label{sufficientcondforlust} Let $X$ be a Banach
  space with a Schauder basis~$(b_n)_{n\in\N_0}$, and let $C\in
  [1,\infty)$. Suppose that $X$ contains a sequence $(F_n)_{n\in\N_0}$
  of finite-dimensional subspaces satisfying
  \begin{equation}\label{sufficientcondforlusteq1}
    b_0,b_1,\ldots,b_n\in F_n\qquad\text{and}\qquad  \uc(F_n)\le C\qquad
    (n\in\N_0).  
  \end{equation}
  Then $X$ has l.u.st.\ with constant at most~$C+\delta$ for
  each~$\delta>0$.
\end{lemma}
\beginpf 
Take $\epsilon\in(0,\frac12)$ such that $C/(1-2\epsilon)<C+\delta$,
and let~$E$ be a $k$-di\-men\-sional sub\-space of~$X$ for some
$k\in\N$. Approximation of each vector of an Auerbach basis for~$E$
shows that, for each $\eta>0$, there is $M\in\N_0$ such that
\begin{equation}\label{sufficientcondforlustEq2}
  \| x - P_mx\|\le\eta\|x\|\qquad (m\ge M,\, x\in E), \end{equation}
where $P_m$ denotes the $m^{\text{th}}$ basis projection associated
with~$(b_n)_{n\in\N_0}$.  
Applying this conclusion with $\eta>0$ chosen such
that  
\begin{equation}\label{sufficientcondforlusteq3}
\frac{\eta\sqrt{k}}{1-\eta}\le\frac{\epsilon}{1-\epsilon}, 
\end{equation}
we obtain by Remark~\ref{remark_bddbelowoperators} that the operator
$U\colon x\mapsto P_Mx,\ E\to P_M(E)$, is an isomorphism with $\|
U\|\le 1 +\eta$ and $\|U^{-1}\|\le (1 -\eta)^{-1}$.  

Since $U(E) = P_M(E)\subseteq\spa\{b_0,b_1,\ldots,b_M\}\subseteq F_M$
and $\dim U(E) = k$, we can find a projection~$Q$ on~$F_M$ such that
$Q(F_M) = U(E)$ and $\| Q\|\le\sqrt{k}$ by the Kadec--Snobar theorem
(\emph{e.g.}, see~\cite[Theorem~4.18]{djt}).  The operator
\mbox{$T\colon x\mapsto x - Qx + U^{-1}Qx,\ F_M\to X$}, then satisfies
\[ \| x - Tx\| = \| Qx - U^{-1}Qx\| = \| P_MU^{-1}Qx - U^{-1}Qx\|\le
\eta\|U^{-1}Qx\|\le \eta\|U^{-1}\|\,\|Q\|\,\|x\|, \] where the
penultimate estimate follows from~\eqref{sufficientcondforlustEq2},
and hence $\| x - Tx\|\le\epsilon(1-\epsilon)^{-1}\|x\|$ for each
$x\in F_M$ by~\eqref{sufficientcondforlusteq3}. Since
$\epsilon(1-\epsilon)^{-1}<1$, Remark~\ref{remark_bddbelowoperators}
implies that~$T$ is an isomorphism onto its image, and
\[ d_{\text{BM}}\bigl(F_M,T(F_M)\bigr)\le\frac{1 +
  \epsilon(1-\epsilon)^{-1}}{1 - \epsilon(1-\epsilon)^{-1}} =
\frac{1}{1-2\epsilon}, \] so $\uc\bigl(T(F_M)\bigr)\le
C/(1-2\epsilon)\le C+\delta$ by~\eqref{remarklusteq1}. 

The conclusion now follows because $E\subseteq T(F_M)$. Indeed, for
each $x\in E$, $y := Ux$ belongs to~$F_M$ and satisfies $Qy = y$, so
that
\[ T(F_M)\ni Ty = y - Qy + U^{-1}y = x, \] as desired. (In fact, an
easy dimension argument shows that $T(F_M) = \ker Q + E$.)  \eopf

\begin{proposition}\label{Jpinftyfailslust}
  The Banach space $J_p^{(\infty)}$ does not have l.u.st.\ for any
  $p\in(1,\infty)$.
\end{proposition}

\beginpf 
Assume towards a contradiction that $J_p^{(\infty)}$ has l.u.st.\ with
constant at most \mbox{$C\ge 1$} for some $p\in(1,\infty)$, and let
$n\in\N_0$. Denote by $\iota_n\colon J_p^{(n)}\to J_p^{(\infty)}$ and
\mbox{$\rho_n\colon J_p^{(\infty)}\to J_p^{(n)}$} the canonical
$n^{\text{th}}$ coordinate embedding and projection, respectively, and
let $j_n\colon J_p^{(n)}\to J_p$ be the natural inclusion operator. By
assumption, $\iota_n(J_p^{(n)})$ is contained in some
finite-dimensional subspace~$F_n$ of~$J_p^{(\infty)}$ with
$\uc(F_n)\le C$.

Let $R_n\colon J_p\to J_p$ be the $(n+2)$-fold right shift given by
$R_ne_k = e_{n+k+2}$ for each $k\in\N_0$. This defines an operator of
norm~$\sqrt[p]{2}$ on~$J_p$, and $R_n$ is bounded below by~$1$.
Lemma~\ref{factsaboutcomplsubspacesoJp}%
\romanref{factsaboutcomplsubspacesoJp3} implies that there are
operators \mbox{$U\in\allop(J_p^{(\infty)},J_p)$} and
$V\in\allop(J_p,J_p^{(\infty)})$ such that \mbox{$VU =
  I_{J_p^{(\infty)}}$}; we may clearly suppose that $V$ has norm one.

We shall now consider the operator
\[ S_n := j_n\rho_n + R_nU(I_{J_p^{(\infty)}} -
\iota_n\rho_n)\in\allop(J_p^{(\infty)},J_p). \]
The obvious norm
estimates show that $\|S_n\|\le 1+\sqrt[p]{2}\,\|U\|$.  To prove that
$S_n$ is bounded below by~$1$, let $x\in J_p^{(\infty)}$ and
$\epsilon>0$ be given. Introducing \mbox{$y := (I_{J_p^{(\infty)}} -
  \iota_n\rho_n)x\in J_p^{(\infty)}$}, we obtain
\begin{equation}\label{sufficientcondforlusteq2}
  \| x\|_{J_p^{(\infty)}}^p = \|\rho_nx\|_{J_p}^p +
\|y\|_{J_p^{(\infty)}}^p = \|j_n\rho_nx\|_{J_p}^p +
\|VUy\|_{J_p^{(\infty)}}^p\le \|j_n\rho_nx\|_{J_p}^p 
  + \|R_nUy\|_{J_p}^p
\end{equation} 
because $\|V\| = 1$ and $R_n$ is bounded below by~$1$.  Since
$j_n\rho_nx\in\spa\{e_0,e_1,\ldots, e_n\}$, there is a subset~$A$
of~\mbox{$[0,n+1]\cap\N$} such that $\|j_n\rho_nx\|_{J_p} =
\nu_p(j_n\rho_nx,A)$.  Similarly, as
$R_nUy\in\closedspa\{e_{n+2},e_{n+3},\ldots\}$, we can find a finite
subset~$B$ of $[n+1,\infty)\cap\N$ such that $\|R_nUy\|^p_{J_p}\le
\nu_p(R_nUy,B)^p + \epsilon$. Combining these identities
with~\eqref{sufficientcondforlusteq2}, we conclude that
\begin{align*} \| x\|_{J_p^{(\infty)}}^p -\epsilon &\le
  \nu_p(j_n\rho_nx,A)^p + \nu_p(R_nUy,B)^p \le \nu_p\bigl(j_n\rho_nx +
  R_nUy,A\cup B\bigr)^p\le\|S_nx\|_{J_p}^p,
\end{align*}  
and letting $\epsilon$ tend to~$0$, we see that $S_n$ is bounded
below by~$1$, as stated.

Thus Remark~\ref{remark_bddbelowoperators} implies that
$d_{\text{BM}}\bigl(F_n,S_n(F_n)\bigr)\le \|S_n\|\le 1+\sqrt[p]{2}\|U\|$,
and therefore $\uc\bigl(S_n(F_n)\bigr)\le
C\bigl(1+\sqrt[p]{2}\|U\|\bigr)$ by~\eqref{remarklusteq1}. Moreover,
for each $k\in\{0,1,\ldots,n\}$, we have $\iota_ne_k\in F_n$, so that
$S_n(F_n)\ni S_n(\iota_ne_k) = e_k$ because $j_n\rho_n\iota_ne_k = e_k$
and $(I_{J_p^{(\infty)}} - \iota_n\rho_n)\iota_n = 0$. 

Hence the sequence $\bigl(S_n(F_n)\bigr)_{n\in\N_0}$ satisfies both
parts of~\eqref{sufficientcondforlusteq1}, so
Lemma~\ref{sufficientcondforlust} implies that~$J_p$ has l.u.st.,
contra\-dict\-ing the above-mentioned theorem of Johnson and Tzafriri
that this is impossible for a quasi-reflexive Banach space. \eopf

\begin{corollary}\label{GpnotisomorphictoJpinfty}
  The Banach spaces $G_p$ and $J_p^{(\infty)}$ are not isomorphic for
  any $p\in(1,\infty)$.
\end{corollary}
\beginpf This is clear because, as remarked above, $G_p$ has an
unconditional basis and thus l.u.st., whereas $J_p^{(\infty)}$ does
not by Proposition~\ref{Jpinftyfailslust}. \eopf \bigskip

The \textsl{proof of Theorem~\ref{newoperatoridealonJp}} is now
easy. Recall that $\weaklycompactop(J_p) =
\closedopidealg_{J_p^{(\infty)}}(J_p)$. The inclusions
$\closedopidealg_{\ell_p}(J_p)\subsetneq \closedopidealg_{G_p}(J_p)$
and $\closedopidealg_{G_p}(J_p)\subsetneq
\closedopidealg_{J_p^{(\infty)}}(J_p)$ both follow from
Lemma~\ref{lemmaopidealsandcomplementedsubspaces}, which applies by
Lemma~\ref{factsaboutcomplsubspacesoJp} and the facts that
$\ell_p\not\cong G_p$ and $G_p\not\cong J_p^{(\infty)}$. The second of
these facts was proved in Corollary~\ref{GpnotisomorphictoJpinfty},
while the first can be justified in various ways; for instance,
$\ell_p$ is uniformly convex with type~$\min\{2,p\}$ and
cotype~$\max\{2,p\}$, whereas $G_p$ is not uniformly convexifiable,
has type~$1$ and fails to have finite cotype. \eopf

\section*{Acknowledgement}
The third author is grateful to Nigel Kalton and Charles Read for
helpful conversations regarding the approach taken in
Section~\ref{section_pfofThmnewoperatoridealonJp} to prove that the
Banach space $J_p^{(\infty)}$ does not have l.u.st.

\bigskip

\begin{center}
  Department of Mathematics and
  Statistics,   Fylde College\\
  Lancaster University,
  Lancaster LA1 4YF, UK;\\[2ex]
  e-mail: \texttt{alistairbird@gmail.com},
  \texttt{g.jameson@lancaster.ac.uk} and
  \texttt{n.laustsen@lancaster.ac.uk}
\end{center}

\end{document}